\documentclass[fleqn,runningheads]{svjour2}

\smartqed

\usepackage{amsmath}
\usepackage{natbib}
\usepackage{amsfonts}
\usepackage{amssymb}

\DeclareMathOperator{\rank}{rank}
\DeclareMathOperator{\ord}{ord}
\DeclareMathOperator{\C}{\mathbb{C}}
\DeclareMathOperator{\charact}{char}
\DeclareMathOperator{\Quot}{Quot}
\DeclareMathOperator{\Char}{char}

\journalname{Applicable Algebra Eng., Comm. and Comp.}

\begin{document}

\sloppy

\title{Canonical Characteristic Sets of Characterizable Differential Ideals \thanks{The work was partially supported by the Russian Foundation for Basic Research, project no. 05-01-00671, by NSF Grant CCR-0096842, and by NSERC Grant PDF-301108-2004.}}

\author{Oleg Golubitsky
\and
Marina Kondratieva
\and
Alexey Ovchinnikov}

\institute{
Oleg Golubitsky \at
University of Western Ontario\\
Department of Computer Science\\
London, Ontario, Canada N6A 5B7\\
\email{oleg.golubitsky@gmail.com}\\
\and
Marina V. Kondratieva \at
Moscow State University\\
Department of Mechanics and Mathematics\\
Leninskie gory, Moscow, Russia, 119992\\
\email{kondra\_m@shade.msu.ru}\\
\and
Alexey Ovchinnikov\thanks{This author was also partially supported by NSF Grant CCR-0096842.}\\
North Carolina State University\\
Department of Mathematics\\
Raleigh, NC 27695-8205, USA\\
\email{aiovchin@ncsu.edu}
}

\date\today

\maketitle

\begin{abstract}
We study the concept of {\it canonical characteristic set\ } of a characterizable differential ideal.
We propose an efficient algorithm that transforms any characteristic set into the canonical one.
We prove the basic properties of canonical characteristic sets. In particular,
we show that in the ordinary case for any ranking the order of each element of the canonical characteristic set of a characterizable
differential ideal is bounded by the order of the ideal.
Finally, we propose a factorization-free algorithm for computing the canonical characteristic set of a characterizable differential
ideal represented as a radical ideal by a set of generators. The algorithm is not restricted to the ordinary case
and is applicable for an arbitrary ranking.
\end{abstract}



\newtheorem{algorithm}{Algorithm}

\newcommand{\Le}{\leqslant}
\newcommand{\Ge}{\geqslant}

\begin{section}{Introduction}
Characterizable ideals have been introduced by Hubert \cite{Fac}
and play a central role in the constructive theory of radical
differential ideals. On the one hand, characterizable ideals
can be specified by their characteristic sets and this representation
allows to solve important computational problems such as ideal
membership. On the other hand, there exist factorization-free
algorithms that decompose any radical differential ideal specified
by a system of generators into an intersection of characterizable ideals
represented by their characteristic sets, which allows to reduce
computational problems for radical ideals to characterizable ones.

Unlike reduced Gr\"obner bases of algebraic ideals, characteristic sets of
(differential) ideals are not unique. The problem of constructing a
unique, or canonical, representation was first addressed in
\cite{Bou2} in the case of regular differential ideals, for which 
the concept of characteristic presentation (which is almost a characteristic set)
is defined. It is shown that, whenever a regular
differential ideal has a characteristic presentation, it is unique.

In \cite{Fac}, the concept of characteristic
presentation is replaced by that of a characterizing set
of a characterizable differential ideal, and an algorithm for
computing such a set is proposed. Further development and
generalization to triangular sets
appeared in \cite{Pol,Dif}.  The efficiency of the mentioned algorithm
was improved in \cite{Hubert04}. 

%One can download the related software 
%from {\tt http://www-sop.inria.fr/cafe/Evelyne.Hubert/diffalg/}.

Finally, in \cite{Bou3}, an additional requirement on characterizing
sets is imposed: the initials of the polynomials in the characteristic 
set are required to be independent of the leaders and the polynomials
themselves to be primitive over the ring of polynomials in
non-leaders. It is shown that, for characterizable ideals, the
characteristic set satisfying this requirement is unique up to 
multiplication by elements of the basic field. In the algebraic
case, a triangular set satisfying the above requirement is called
a reduced Gr\"obner chain \cite{Pol}, and it shown in
\cite[Proposition 5.17]{Pol} 
that every algebraic characterizable ideal admits a unique characteristic
set that is a reduced Gr\"obner chain. By applying \cite[Theorem
  5.5]{Dif}, this result can be trivially lifted to differential 
characterizable ideals. 

We give afresh a detailed constructive proof of existence and
uniqueness of the canonical characteristic set of a characterizable
differential ideal. Based on the algorithm for
inverting the polynomials w.r.t. a characteristic set
\cite{Unmixed}, we propose an algorithm
(Algorithm~\ref{CanonicalAlgorithm}) that, given any characteristic set
that characterizes an ideal, constructs the canonical characterizing
set. These results are contained in Sections~\ref{Definitions} and 
\ref{CanonicalProperties}.

The canonical characteristic set is a convenient tool for testing
equality of characterizable differential ideals. Indeed, two
characterizable ideals coincide iff their canonical characteristic
sets coincide. This method is used in \cite{ComputingResolvent}
in order to make a certain characteristic decomposition irredundant.
On the other hand, equality of prime differential
ideals can be tested using the method described in \cite[Exercise 1,
  page 171]{Kol}. As we show in Section~\ref{Equality}, this method 
is also applicable to characterizable ideals.

We study the essential properties of canonical characteristic sets 
(see Propositions~\ref{Irreducible}, \ref{LowestRank}, and
\ref{Separants}). In particular, we prove that the partial derivatives 
of the polynomials in the canonical characteristic sets w.r.t. any 
differential indeterminate do not belong to the ideal. Although these 
properties look algebraic, we also bring some differential flavor into 
their study. We also describe the structure of all characteristic sets 
whose initials do not depend on leaders
(Corollary~\ref{Multiples}).

%In order to obtain additional differential properties of canonical 
%characteristic sets, we first study what properties one can guarantee
%for some characteristic sets of characterizable ideals. 
%Section~\ref{Result} is essentially about this. 

The main results of the paper are the following.
First, we show that  
{\it in the ordinary case for any ranking and a prime differential ideal there exists a 
characteristic set with bounded in advance orders of its elements} 
(Theorem~\ref{OrderlyElimination}).  Our result generalizes \cite[Theorem 24]{Sadik}, which was proved there for {\it elimination rankings only}.
Then we obtain a similar result 
for {\it characterizable differential ideals} 
(see Theorem~\ref{CharacterizableBound}) without any restrictions to rankings. In 
Section~\ref{CanonicalProperties} we show how to apply this result. 
In Theorem~\ref{HubertCanonical} we obtain a {\it bound on the orders of 
the elements in the canonical characteristic set} of a characterizable 
differential ideal.

Considering prime differential ideals as a particular case of characterizable
ones, in Example~\ref{Counterexample} we show that Sadik's property of 
irreducibility \cite{Sadik} does not
necessarily hold for the canonical characteristic set. This property
is very important for the proof of Theorem~\ref{OrderlyElimination}
together with another natural property giving a characteristic set of
a prime differential ideal
with a predicted bound for the order of its elements.

As it is demonstrated in \cite{Bou4},
it is possible to convert characteristic sets of prime differential
ideals from one ranking to another efficiently. Furthermore, an
efficient Monte-Carlo algorithm for converting characteristic sets
of prime algebraic ideals has been developed in \cite{DahanPoster}. 
This algorithm can be applied in case of prime ordinary
differential ideals using the reduction of this case 
to the algebraic one carried out in \cite{Bounds}; the reduction
essentially involves the concept of canonical characteristic set and its
properties. 

Finally, we propose an algorithm for computing the canonical characteristic set
of a characterizable differential ideal from a set of its generators (Algorithm~\ref{CharCharComputation}). Examples
in Section~\ref{Examples} illustrate our algorithms, as well as justify why the requirement of characterizability
cannot be relaxed to that of radicality.
\end{section}

\begin{section}{Preliminaries}\label{Prelim}

\begin{subsection}{Basic definitions}
Differential algebra studies systems of polynomial partial differential equations from
the algebraic point of view. The approach is based on the concept of differential ring introduced by Ritt.
Recent tutorials on the constructive theory of differential ideals are presented in \cite{Dif, Sit}. We also use the Gr\"obner bases technique \cite{Bb,Wei}.
A differential ring is a commutative ring with unity endowed
with a set of differentiations $\Delta = \{\delta_1,\ldots,\delta_n\}$. The case of $\Delta = \{\delta\}$ is called {\it ordinary}. If $R$ is an ordinary differential ring and $y \in R$, we denote $\delta^ky$ by $y^{(k)}$.

 Construct the multiplicative monoid $\Theta = (\delta_1^{k_1}\delta_2^{k_2}\cdots\delta_n^{k_n}, k_i \Ge 0).$ The ring of {\it differential polynomials} in $l$ differential indeterminates over a differential field $k$ is a ring of commutative polynomials with coefficients in $k$ in the infinite set of variables $\{\Theta y_i$, $1\Le i \Le l\}$ (see \cite{Kol, Pan, Rit}). Derivatives from $\Theta$ act on variables as $\theta_1(\theta_2y_i) = (\theta_1\theta_2)y_i$ for all $\theta_1,\theta_2 \in \Theta$ and $1 \Le i \Le l$. The ring of differential polynomials is denoted by $k\{y_1,\dots,y_l\}$ or $k\{Y\}$, where $Y=\{y_1,\dots,y_l\}$ is the set of differential indeterminates. We consider the case of $\charact k=0$ only. We denote polynomials by $f, g, h, \ldots$ and use letters $I, J, P$ for ideals.

Let $F \subset k\{y_1,\dots,y_l\}$ be a set of differential polynomials. For the differential and radical differential ideal generated by $F$ in $k\{y_1,\dots,y_l\}$, we use notations $[F]$ and $\{F\}$, respectively. A prime differential ideal containing radical differential
ideal $\{F\}$ is called
a prime component of $\{F\}$. A prime component is called {\it essential}, if it is not contained in any other prime component of $\{F\}$.
According to \cite[Section I.16]{Rit}, every radical differential ideal has finitely many essential prime components and is equal to
their intersection. Moreover, any finite prime decomposition of $\{F\}$ contains all its essential components.

We need the notion of reduction for algorithmic computations. First, we introduce a {\it ranking} on the set of differential variables of $k\{y_1,\ldots,y_l\}$.  A ranking is a total ordering on the set $\{\theta y_i\}$, where $\theta \in \Theta$ and $1\Le i \Le l$,
satisfying the following conditions:
\begin{enumerate}
\item $\theta u \Ge u,$
\item $u \Ge v \Longrightarrow \theta u \Ge \theta v.$
\end{enumerate}

Let $u$ be a differential variable in $k\{y_1,\ldots,y_l\}$, that is, $u = \theta y_j$ for a differential operator $\theta = \delta_1^{k_1}\delta_2^{k_2}\cdots\delta_n^{k_n} \in \Theta$ and $1\Le j \Le l$. The {\it order} of $u$ is defined as $\ord u=\ord\theta=k_1+\ldots+k_n$. If $f$ is a differential polynomial then $\ord f$ denotes the maximal order of differential variables appearing effectively in $f$. A ranking $>$ is said to be {\it orderly} iff $\ord u > \ord v$ implies $u > v$ for all differential variables $u$ and $v$. A ranking $>_{el}$ is called an {\it elimination} ranking iff $y_i >_{el} y_j$ implies $\theta_1y_i >_{el} \theta_2y_j$ for all $\theta_1, \theta_2 \in \Theta$.

Let a ranking $<$ be fixed.
The differential variable $\theta y_j$ of the highest rank appearing in a differential polynomial $f \in k\{y_1,\dots,y_l\} \setminus k$ is called the {\it leader} of $f$. We denote the leader by $u_f$. Represent $f$ as a univariate polynomial in $u_f$:
$$
f = I_f u_f^n + a_1 u_f^{n-1} + \ldots + a_n.
$$
The polynomial $I_f$ is called the {\it initial} of $f$.

Apply any $\delta \in \Delta$ to $f$:
$$ \delta f = \frac{\partial f}{\partial u_f}\delta u_f + \delta I_f u_f^n +
\delta a_1 u_f^{n-1}+\ldots + \delta a_n.
$$
The leading variable of $\delta f$ is $\delta u_f$ and the initial of $\delta f$ is called the {\it separant} of $f$, denoted $S_f$. Note that for all $\theta \in \Theta,$ $\theta \ne 1,$ the initial of $\theta f$ is equal to $S_f$. The differential monomial $u_f^n$ is called the {\it rank} of $f$ and denoted $\rank f$.

Extend the ranking relation on differential variables to ranks: $u_1^{d_1}>u_2^{d_2}$ iff either $u_1>u_2$ or
$u_1=u_2$ and $d_1>d_2$.
Also, a ranking $>$ on differential variables induces a lexicographic ordering $>_{\rm lex}$ on their power products.
This ordering can be extended to a partial ordering on differential polynomials, which we denote by $>_{\rm lex}$ as well.

We say that a differential polynomial $f$ is {\it partially reduced} w.r.t. $g$ iff no proper derivative of $u_g$ appears in $f$. A differential polynomial $f$ is {\it reduced} w.r.t. a differential polynomial $g$ iff $f$ is partially reduced w.r.t. $g$ and $\deg_{u_g} f < \deg_{u_g} g$. Consider any subset $\mathbb{A} \subset k\{y_1,\ldots,y_l\}\setminus k$. We say that $\mathbb{A}$ is autoreduced iff each element of $\mathbb{A}$ is reduced w.r.t. all the others. Every autoreduced set is finite \cite[Chapter I, Section 9]{Kol}. For autoreduced sets we use capital letters $\mathbb{A, B, C,}$ \ldots and notation $\mathbb{A} = A_1,\ldots,A_p$ to specify the list of the elements of $\mathbb{A}$ arranged in order of increasing rank.

We denote the product of the initials and the separants of the elements of $\mathbb{A}$ by $I_\mathbb{A}$ and $S_\mathbb{A}$, respectively. Denote $I_\mathbb{A}\cdot S_\mathbb{A}$ by $H_\mathbb{A}$. Let $S$ be a finite set of differential polynomials. Denote by $S^\infty$ the multiplicative set containing $1$ and generated by $S$. Let $I$ be an ideal in a commutative ring $R$. The {\it saturated ideal} $I:S^\infty$ is defined as $\{a \in R\:|\:\exists s \in S^\infty: sa \in I\}$. If $I$ is a differential ideal then $I:S^\infty$ is also a differential ideal (see \cite{ Kol, Rit, Pan, Sit}).

Consider two differential polynomials $f$ and $g$ in $R = k\{y_1,\ldots, y_l\}$. Let $I$ be the differential ideal in $R$ generated by $g$. Applying a finite number of differentiations and pseudo-divisions (multiplying by initials and separants together with differentiations and algebraic reductions) one can compute a {\it differential partial remainder} $f_1$ and a {\it differential remainder} $f_2$
of $f$ w.r.t. $g$ such that there exist $s \in S_g^\infty$ and $h \in H_g^\infty$ satisfying $sf \equiv f_1$ and $hf \equiv f_2 \mod I$ with $f_1$ and $f_2$ partially reduced and reduced w.r.t. $g$, respectively (see \cite{Fac} for definitions and the algorithm for computing remainders).

Let $\mathbb{A}$ be an autoreduced set in $k\{y_1,\ldots,y_l\}$. Consider the polynomial ring $k[x_1,\ldots,x_n]$ with $x_1,\ldots,x_n$ belonging to $\Theta Y$ for $Y = y_1,\ldots,y_l$. Let $L, N \subset \{x_1,\ldots,x_n\}$ be the sets of ``leaders'' and ``non-leaders'' appearing in the autoreduced set $\mathbb{A}$, respectively. We denote $k[x_1,\ldots, x_n]$ by $k[N][L]$ and the leader of $A_i$ by $u_{A_i}$ or $u_i$ for each $1 \Le i \Le p$. In what follows,
we will often consider elements of $\mathbb{A}$ as polynomials in leaders with coefficients being polynomials
in non-leaders.

Let $\mathbb{A} = A_1,\ldots,A_r$ and $\mathbb{B} = B_1,\ldots,B_s$ be autoreduced sets. We say that $\mathbb{A}$ has lower rank than $\mathbb{B}$ iff there exists $k \Le r, s$ such that $\rank A_i$ = $\rank B_i$ for $1 \Le i < k$ and $\rank A_k < \rank B_k$, or if $r > s$ and $\rank A_i = \rank B_i$ for $1 \Le i \Le s$. We say that $\rank\mathbb{A} = \rank\mathbb{B}$ iff $r=s$ and $\rank A_i = \rank B_i$ for $1 \Le i \Le r$.

The following notion of a characteristic set in {\it Kolchin's sense} in characteristic zero is crucial in our further discussions. It was first introduced by Ritt for prime differential ideals, and then extended by Kolchin to arbitrary differential ideals.

\begin{definition} \cite[page 82]{Kol} An autoreduced set of the lowest rank in an ideal $I$ is
called a {\it characteristic set} of $I$.
\end{definition}

We call these sets {\it Kolchin characteristic sets} to avoid confusion with other notions, e.g., in \cite{Fac,Dif} characteristic sets are used in Kolchin's sense and in some other senses.

As it is mentioned in \cite[Lemma 8, page 82]{Kol}, in the case of $\Char k = 0$, a set $\mathbb{A}$ is a characteristic set of a proper differential ideal $I$ iff each element of $I$ reduces to zero w.r.t. $\mathbb{A}$. Moreover, the leaders and the correspondent degrees of these leaders of any two characteristic sets of $I$ coincide.

Let $v$ be a derivative in $k\{y_1,\ldots,y_l\}.$ Denote by $\mathbb{A}_v$ the set of the elements of $\mathbb{A}$ and their derivatives of rank strictly lower than $v$.
\begin{definition}\cite[III.8]{Kol}
An autoreduced set $\mathbb{A}$ is called {\it coherent} if whenever $A, B \in \mathbb{A}$ are such that $u_{A}$ and $u_{B}$ have a common derivative $v = \psi u_{A} = \phi u_{B}$, then $S_{B}\psi A - S_{A}\phi B \in (\mathbb{A}_v):H_\mathbb{A}^\infty.$
\end{definition}

Any characteristic set of a differential ideal is coherent (see \cite{Kol,Rit, Pan, Sit}).

\begin{definition}\cite[Definition 2.6]{Fac} A differential ideal $I$ in $k\{y_1,\ldots,y_l\}$ is
said to be {\it characterizable} if there exists a characteristic set $\mathbb{A}$ of $I$ in Kolchin's sense
such that $I = [\mathbb{A}]:H_\mathbb{A}^\infty.$ We call any such characteristic set $\mathbb{A}$ a {\it characterizing} set
of $I$.
\end{definition}

Characterizable ideals are radical \cite[Theorem 4.4]{Fac}.

Let $\mathbb{A}$ be an autoreduced set in $k\{y_1,\ldots,y_l\} = k\{Y\}$, and let $k[N][L]$ be the polynomial ring associated with $\mathbb{A}$, where $L$ is the set of
leaders of polynomials in $\mathbb{A}$ and $N$ is the set of
non-leaders, i.e., $N = \Theta Y\setminus \Theta L.$ Note that the set $N$ can be infinite when $\Delta \ne \varnothing$.

\begin{definition}\label{D1}\cite[Definition 1.2.1]{Unmixed}
Let $f$ be a polynomial in
$k[N][L]$. We say that $f$
is {\it invertible} w.r.t. $\mathbb{A}$
iff $f$ is invertible modulo $(\mathbb{A})$ in $k(N)[L]$, i.e., there exist $g
\in
k[N,L]$ and $h \ne 0$ in $k[N]$
such that $f\cdot g \equiv h \mod (\mathbb{A})$.
\end{definition}

We say that the initials of an autoreduced set $\mathbb{A}=A_1,\ldots,A_p$ are {\it invertible} if $I_{A_i}$ is invertible w.r.t. the autoreduced set $A_1,\ldots,A_{i-1}$ for all $i$, $2\Le i \Le p$. We also say that the separants of $\mathbb{A}$ are {\it invertible} if $S_{A_i}$ is invertible w.r.t. the autoreduced set $A_1,\ldots,A_{i}$ for all $i$, $1\Le i \Le p$.

\end{subsection}

\begin{subsection}{Important assertions}
Consider several important results concerning radical differential ideals in rings of differential polynomials. The technique described in \cite{Fac, Kol} helps us to cover some properties of these ideals.

\begin{theorem}\cite[III.8, Lemma 5]{Kol}\label{Rosenfeld} Let $\mathbb{A}$ be a coherent autoreduced set in $k\{y_1,\ldots,y_l\}$. Suppose that a differential polynomial $g$ is partially reduced w.r.t. $\mathbb{A}$. Then $g \in [\mathbb{A}]:H_\mathbb{A}^\infty$ iff $g \in (\mathbb{A}):H_\mathbb{A}^\infty$.
\end{theorem}

Note that Theorem~\ref{Rosenfeld} is also known as {\it Rosenfeld's lemma}.

\begin{theorem}\cite[Theorem 3.2]{Fac}\label{T20} Let $\mathbb{A}$ be an autoreduced set of $k[N][L]$. If $1 \notin (\mathbb{A}):H_\mathbb{A}^\infty$ then any minimal prime of $(\mathbb{A}):H_\mathbb{A}^\infty$ admits the set of non-leaders of
$\mathbb{A}$, $N$, as a transcendence basis. More specially, any characteristic set of a minimal prime of $(\mathbb{A}):H_\mathbb{A}^\infty$ has the same set of leaders as $\mathbb{A}$.
\end{theorem}

\begin{theorem}\cite[Theorem 4.5]{Fac}\label{Prime} Let $\mathbb{A}$ be a coherent autoreduced set of $R = k\{y_1,\ldots,y_l\}$ such that $1 \notin [\mathbb{A}]:H_\mathbb{A}^\infty$. There is a one-to-one correspondence between the minimal primes of $(\mathbb{A}):H_\mathbb{A}^\infty$ in $k[N][L]$ and the essential prime components of $[\mathbb{A}]:H_\mathbb{A}^\infty$ in
$R$. Assume $\mathbb{C}_i$ is a characteristic set of a minimal prime of $(\mathbb{A}):H_\mathbb{A}^\infty$. Then $\mathbb{C}_i$ is the characteristic set of a single essential prime component of $[\mathbb{A}]:H_\mathbb{A}^\infty$ (and vice versa).
\end{theorem}

\begin{theorem}\cite[Theorem 1.2.2]{Unmixed}\label{KandriInvertible}
An autoreduced set $\mathbb{A} = A_1,\ldots,A_p$ is a characteristic set of the algebraic ideal $(\mathbb{A}):I_\mathbb{A}^\infty$ in $k[N,L]$ iff the initials $I_{A_i}$ are invertible for all $2 \Le i \Le p$.
\end{theorem}

\begin{theorem}\cite[Theorem 2.1.1]{Unmixed}\label{KandriInvertibleSep}
A coherent autoreduced set $\mathbb{A} = A_1,\ldots,A_p$ is a characteristic set of the  ideal $[\mathbb{A}]:H_\mathbb{A}^\infty$ iff the initials $I_{A_i}$ and $S_{A_j}$ are invertible for all $2 \Le i \Le p,$ $1 \Le j \Le p$.
\end{theorem}

\begin{lemma}\label{LoseVar} Let $\mathbb{A} = A_1, \ldots, A_p$ be an autoreduced set in the ring $k[x_1,\ldots,x_m] = R$ and a characteristic set of $(\mathbb{A}):I_\mathbb{A}^\infty$. Suppose that a polynomial $f = a_mx_t^m + \ldots + a_0 \in R$ is reducible to zero w.r.t. $\mathbb{A}$ and the indeterminate $x_t$ does not appear in $A_i$ for each $1\Le i \Le p$. Then $a_j$ is reducible to zero w.r.t. $\mathbb{A}$ for all $0 \Le j \Le m$.
\end{lemma}
\begin{proof}
Since $f$ is reducible to zero w.r.t. $\mathbb{A}$, there exists $I \in  I_\mathbb{A}^\infty$ such that
$$
I\cdot f = \sum\limits_{i=1}^p g_iA_i.
$$

Let $g_i = \sum\limits_{j=1}^{t_i}h_{i_j} x_t^j$ for each $1 \Le i \Le p$. Thus, we have $I\cdot \sum\limits_{k=0}^m a_k x_t^k = \sum\limits_{k=0}^q d_k x_t^k$ with $d_k \in (\mathbb{A})$. Hence, $I\cdot a_i \in (\mathbb{A})$ for each $1 \Le i \Le m$,
that is, $a_i \in (\mathbb{A}):I_\mathbb{A}^\infty$. Since  $\mathbb{A}$ is a characteristic set of $(\mathbb{A}):I_\mathbb{A}^\infty$, we have that all $a_i$ are reducible to zero w.r.t. $\mathbb{A}$.
\end{proof}

We note an important property of a characteristic set
of an arbitrary differential ideal.

\begin{proposition}\label{Multiplication} Let $I$ be a differential ideal with two characteristic sets $\mathbb{A} = A_1,\ldots,A_p$ and $\mathbb{C} = C_1,\ldots,C_p$. If the initials of both $\mathbb{A}$ and $\mathbb{C}$ do not depend on the leaders then
$$I_{A_i}C_i = I_{C_i}A_i$$
for all $i,$ $1\Le i \Le p$.
\end{proposition}
\begin{proof}
Let
\begin{align*}
A_i &= I_{A_i}u_{A_i}^{n_i} + a_{i,1}u_{A_i}^{n_i-1} + \ldots + a_{i,n_i}, \\
C_i &= I_{C_i}u_{A_i}^{n_i} + c_{i,1}u_{A_i}^{n_i-1} + \ldots + c_{i,n_i},
\end{align*}
where $u_{A_i}$
is the leader of both $A_i$ and $C_i$, $1\Le i \Le p$.
We have
$$ f_i := I_{A_i}C_i - I_{C_i}A_i = \sum_{j=1}^{n_i}{(c_{i,j}I_{A_i} -
a_{i,j}I_{C_i})u_{A_i}^{n_i-j}} \in I.$$
By Theorem~\ref{KandriInvertible} we have $\mathbb{A}$ is a characteristic set of the ideal $I$. Hence, the polynomial $f_i$
is reducible to zero w.r.t. $\mathbb{A}$. The initials of $\mathbb{A}$ do not depend on leaders, hence $\mathbb{A}$ is a characteristic set of $(\mathbb{A}):I_\mathbb{A}^\infty$ by Theorem~\ref{KandriInvertible}.
The variable $u_{A_i}$ does not appear in $A_1,\ldots,A_{i-1}$.

Hence, by Lemma~\ref{LoseVar} the coefficients of $f_i$ w.r.t. $u_{A_i}$ are reducible to zero w.r.t. $\mathbb{A}$.
But the initials of both $\mathbb{A}$ and $\mathbb{C}$ do not depend on the leaders,
hence the powers of leaders in the coefficients of $f_i$ are less than or equal to the powers of the corresponding leaders of $A_1,\ldots,A_{i-1}$ and no reduction can go. So, we obtain that
$c_{i,j}I_{A_i} = a_{i,j}I_{C_i},$ $1\Le j\Le n_i$, $1 \Le i \Le p.$
Thus, $f_i =0$ for all $i,$ $1\Le i \Le p.$
\end{proof}

\end{subsection}

\end{section}

\begin{section}{Canonical characteristic sets}\label{Definitions}
\begin{subsection}{Definition and computation}\label{DifferentDefinitions}

Let a differential ranking be {\it fixed}. The following definition is a summary of \cite[Section 2.2.6]{Bou3}.

\begin{definition}\label{SecondDefinition} A characteristic set
  $\mathbb{C} = C_1,\ldots,C_p$ of a characterizable differential 
  ideal $I$ is called {\it canonical} if the following conditions are satisfied:
\begin{enumerate}
\item the initial $I_{C_i}$ depends only on non-leaders $N$ of $\mathbb{C}$ for all $i,$ $1\Le i \Le p$;
\item for each $i$, $1 \Le i \Le p,$
\begin{itemize}
  \item the polynomial $C_i$ has no factors in $k[N]$;
  \item the leading coefficient of the leading monomial of $I_{C_i}$
  w.r.t. the induced lexicographic ordering $>_{\rm lex}$ on monomials
  from $k[N]$ is equal to $1$.
\end{itemize}
\end{enumerate}
\end{definition}

Note that, unlike \cite{Bou3}, we do not require in the definition
that $\mathbb{C}$ characterizes $I$, but we will show
in Proposition~\ref{Irreducible} that this is the case. We have also
replaced the set $N_{\mathbb C}$ of non-leaders effectively occurring
in $\mathbb C$ by the set $N=\Theta Y\setminus \Theta L$ of all non-leaders
(where $L$ is the set of leaders of $\mathbb C$). Clearly, this
replacement yields an equivalent definition, which is more convenient
for us, because it provides a ring $k(N)[L]$ independently
of the choice of the characteristic set $\mathbb C$, while the
field of coefficients $k(N_{\mathbb C})$ of the polynomial ring 
$k(N_{\mathbb C})[L]$ depends on $\mathbb C$.\footnote{The 
idea of constructing a canonical field of coefficients by considering 
the infinite set of all non-leading derivatives was communicated to 
the first author by E. Hubert.}
Other than this detail, the construction of the  
canonical characteristic set follows from \cite[Section 5.4]{Pol}
and \cite[Theorem 5.5]{Dif}.

The following algorithm is a combination of \cite[Algorithm 3.8]{Fac} and
\cite[Algorithm 7.1]{Fac} restricted to characterizable ideals. Its
correctness is justified in Proposition~\ref{Irreducible} and implies
the existence of the canonical characteristic set. The uniqueness of
the canonical characteristic set essentially follows from
\cite[Theorem 3]{Bou3}, yet, since we have
slightly changed the definition,  
we provide another proof of uniqueness. 

\begin{figure}
\begin{algorithm}\label{Canonical}{\sf Canonical Characteristic Set}

{\sc Input:} a characterizing set $\mathbb{A}$ of a characterizable
differential ideal $I = [\mathbb{A}]:H_\mathbb{A}^\infty$ in
$k\{y_1,\ldots,y_l\} = k\{Y\}.$

{\sc Output:} the canonical characteristic set of $I$.

\begin{itemize}
\item $L := $ {\sf Leaders}$(\mathbb{A}).$
\item $N := \Theta Y\setminus \Theta L.$
\item $GB := $ {\sf Reduced\_Gr\"obner\_Basis}$((\mathbb{A}):H_{\mathbb{A}}^\infty)$ in $k(N)[L]$ w.r.t. $>_{\rm lex}$.
%\item $N' := \{x\in N\:|\: x\ \mathrm{appears\ in}\ GB\}.$
\item $\mathbb{C} := $ {\sf Clear\_out\_denominators}$(GB)$ in $k(N)[L].$
\item divide each element of $\mathbb{C}$ by its leading coefficient from $k$
\item Return $\mathbb{C}$.
\end{itemize}
\end{algorithm}
\end{figure}

Note that in the above algorithm the Gr\"obner basis is computed for
a system of polynomials over the field of fractions $k(N)$, where
$N$ is an infinite set. This does not raise any problems
with computability, since the input consists of polynomials
over $k(N_{\mathbb{A}})$, where $N_{\mathbb{A}}$ is a finite subset
of non-leading derivatives effectively present in $\mathbb{A}$. 
Since the Gr\"obner basis computation does not lead out
of the field of definition of the input polynomials, one
can say that, for the the given system $\mathbb{A}$, the 
computation of the Gr\"obner basis of the ideal
$(\mathbb{A}):H_{\mathbb{A}}^\infty$ in $k(N)[L]$ is equivalent
to that in $k(N_{\mathbb{A}})[L]$. 
The same applies to clearing out the denominators 
in the Gr\"obner basis in $k(N)[L]$, i.e., effectively this procedure 
is performed in $k(N_\mathbb{A})[L]$. 

\begin{proposition}\label{Irreducible} Algorithm~\ref{Canonical}
 computes a canonical characteristic set of the characterizable 
 differential ideal $I$. Moreover, this characteristic set 
 characterizes $I$.
\end{proposition}
\begin{proof}
Since ideal $(\mathbb{A}):H_{\mathbb{A}}^\infty$ is characterizable
and according to \cite[Remark after Lemma 3.9]{Fac}, the reduced Gr\"obner basis $GB$ has
$|L|$ elements and characterizes $(\mathbb{A}):H_{\mathbb{A}}^\infty$
in $k(N)[L]$. Thus, $(\mathbb{A}):H_{\mathbb{A}}^\infty=(GB)$
is a zero-dimensional irredundant characteristic decomposition of 
$(\mathbb{A}):H_{\mathbb{A}}^\infty$ in $k(N)[L]$ consisting of the single
component $(GB)$. Now, according to \cite[Theorem 3.10]{Fac},
$(\mathbb{C}):I_{\mathbb{C}}^\infty$ is a single-component characteristic
irredundant decomposition of $(\mathbb{A}):H_{\mathbb{A}}^\infty$ in
$k[N][L]$, and by \cite[Theorem 6.2]{Fac}
$$I= [\mathbb{A}]:H_{\mathbb{A}}^\infty=[\mathbb{C}]:H_{\mathbb{C}}^\infty$$
is a single-component characteristic irredundant decomposition of $I$
in $k\{Y\}$. The latter, in particular, implies that $\mathbb{C}$ 
characterizes $I$.

To clear out denominators means to multiply each $C_i\in GB$ by the
least common multiple $q$ of the denominators of its coefficients. 
Here each
$C_i=u^t+\frac{\alpha_{t-1}}{\beta_{t-1}}u^{t-1}+\ldots+\frac{\alpha_0}{\beta_0}$ 
is considered as a univariate polynomial in its leader $u$, and its
coefficients $\frac{\alpha_i}{\beta_i}$, $i=0,\ldots,t-1$ are assumed
to be irreducible fractions over $k[N]$.

Now let $C_i'=qC_i$. The coefficients of $C_i'$ are polynomials 
$q,\alpha_{t-1}\cdot\frac q{\beta_{t-1}},\ldots,\alpha_0\cdot\frac
q{\beta_0}$, whose greatest common divisor is 1. Indeed, let
$\gamma_i=\frac q{\beta_i}$, and suppose that there exists an
irreducible polynomial $d\in k[N]\setminus k$ such that
$d|q$ and $d|\alpha_i\gamma_i$, $i=0,\ldots,t-1$. If for some 
$j$, $d\not|\gamma_j$, then we have $d|\alpha_j$ and $d|\beta_j=\frac q{\gamma_j}$, which
contradicts the irreducibility of the fraction
$\frac{\alpha_j}{\beta_j}$. Hence, for all $j\in\{0,\ldots,t-1\}$,
$d|\gamma_j$, which implies that
$\beta_j=\frac q{\gamma_j}|\frac qd$. The latter contradicts the fact 
that $q=\mathop{\rm lcm}(\beta_0,\ldots,\beta_{t-1})$.

Thus, set $\mathbb{C}$ satisfies the requirements of
Definition~\ref{SecondDefinition} 
(the leading coefficient in $k$ can be easily canceled) and
hence is the canonical characteristic set of $I$.
\end{proof}

\begin{proposition}
  If $\mathbb C$, $\mathbb D$ are two characteristic sets of a 
  characterizable differential ideal $I$ satisfying the requirements
  of Definition~\ref{SecondDefinition}, then ${\mathbb C}={\mathbb D}$.
\end{proposition}
\begin{proof}
  Let $I_1=({\mathbb C}):H_{\mathbb C}^\infty$ and 
  $I_2=({\mathbb D}):H_{\mathbb D}^\infty$ be the saturated algebraic
  ideals specified by $\mathbb C$ and $\mathbb D$, respectively, 
  considered in the ring $k[L\cup N]$. 

  By the Rosenfeld Lemma, ${\mathbb C}\subset I_2$ and 
  ${\mathbb D}\subset I_1$. 

  Now consider the corresponding ideals $J_1=({\mathbb C}):H_{\mathbb C}^\infty$ and 
  $J_2=({\mathbb D}):H_{\mathbb D}^\infty$ in $k(N)[L]$. 
  Since the polynomials in $k[L\cup N]$ can be also considered as
  elements of $k[N][L]$ and, hence, of $k(N)[L]$, we have
  ${\mathbb C}\subset J_2$ and ${\mathbb D}\subset J_1$. 
  Since $I_{\mathbb C}$ and $I_{\mathbb D}$ belong to $k(N)$, ideals $J_1$
  and $J_2$ are generated by 
  $\mathbb C$ and $\mathbb D$, which implies the equality of these
  ideals. Denote $J=J_1=J_2$.

  Let $\bar{\mathbb C}=\{f/I_f\;|\;f\in{\mathbb C}\}$ and,
  similarly, define $\bar{\mathbb D}$. 
  According to 
  Definition~\ref{SecondDefinition}, the initials of $\mathbb C$
  and $\mathbb D$ depend only on non-leaders, hence $\bar{\mathbb C}$
  and $\bar{\mathbb D}$ are sets of monic polynomials in $k(N)[L]$
  generating the ideal $J$. 
  Moreover, their sets of leading monomials w.r.t. the lexicographic
  ordering $<_{\rm lex}$ on the monomials over $L$ induced by the 
  ranking are equal to $R=\rank{\mathbb C}=\rank{\mathbb D}$. 
  Thus, $\bar{\mathbb C}$ and $\bar{\mathbb D}$ are reduced Gr\"obner 
  bases of $J$, whence they must be equal. 

  Now the conditions of Definition~\ref{SecondDefinition} imply
  that ${\mathbb C}={\mathbb D}$ as well. 
\end{proof}

Note that the equality of ideals $I_1$ and $I_2$, and hence of the
ideals $J_1$ and $J_2$, in the above proof
also follows from \cite[Theorem 6.1]{Bou2}. However, as we have seen,
this equality becomes rather straightforward, once the polynomial ring 
$k(N)[L]$ is considered. 
\end{subsection}

\begin{subsection}{Basic properties}

%Denote $\tilde N=\Theta Y\setminus\Theta L$.

\begin{corollary}\label{Multiples} For a characterizable differential ideal, all characteristic sets with initials containing only non-leaders can be obtained from the canonical one through multiplying its elements by some polynomials from $k[N]$.
\end{corollary}
\begin{proof}
  Let $\mathbb{A}=A_1,\ldots,A_p$ be a characteristic set whose initials belong to $k[N]$. Consider set
  $\mathbb{A}'=A_1',\ldots,A_p'$, where $A_i'=\frac{A_i}{I_{A_i}}\in k(N)[L]$. Then the leading monomials of $\mathbb{A}'$ coincide with
  the ranks of $\mathbb{A}$, which, in turn, coincide with the ranks of the canonical characteristic set $\mathbb{C}$ and with the leading monomials of
  the reduced Gr\"obner basis $GB$ computed in Algorithm~\ref{Canonical}. Moreover, $\mathbb{A}'$ is algebraically autoreduced, hence $\mathbb{A}'=GB$.
  This implies that characteristic set $\mathbb{A}$ can be obtained from the canonical characteristic set $\mathbb{C}$
  through multiplying its elements by polynomials from $k[N]$.

 % Vice versa, if we multiply elements of $\mathbb{C}$ by arbitrary polynomials from %$k[\tilde N]$,
%  $\mathbb{C}$ remains to be autoreduced and preserves its
%  rank, hence it remains to be a characteristic set of the characterizable ideal.
\end{proof}

The next two propositions demonstrate that canonical characteristic sets are ``minimal'' in a certain sense.

\begin{proposition}\label{LowestRank}
Let $\mathbb{C} = C_1,\ldots,C_p$ be the canonical characteristic set of a characterizable differential ideal $I$. Let $\mathbb{B}=B_1,\ldots,B_p$ be any characteristic set  of the ideal $I$ such that the initials of $\mathbb{B}$ depend only on the non-leaders $N$. Then,
$$I_{C_i} \Le_{\rm lex} I_{B_i}$$ for all $i$, $1\Le i \Le p$.
\end{proposition}
\begin{proof} By Corollary~\ref{Multiples} all the other characteristic sets $\mathbb{B}$ with the initials from $k[N]$ can
be obtained from $\mathbb{C}$ by multiplying its elements by polynomials from $k[N]$.
The result follows immediately from this.
\end{proof}

\begin{proposition}\label{Separants} Let $\mathbb{C} = C_1,\ldots,C_p$ be the canonical characteristic set of a characterizable differential ideal $I$. Let $v$ be a differential variable appearing in some $C_i$, $1\Le i \Le p$. Then, $$\dfrac{\partial{C_i}}{\partial{v}} \notin I.$$
\end{proposition}
\begin{proof}
Suppose that $\dfrac{\partial{C_i}}{\partial{v}} \in I.$ Then $v$ appears effectively in the initial $I_{C_i}$.

Indeed, suppose that $v$ is not in $I_{C_i}$, then $\dfrac{\partial{C_i}}{\partial{v}}$ is not reducible w.r.t. $\mathbb{C}$. This contradicts
the fact that $\mathbb{C}$ is a characteristic set of $I$ and $\dfrac{\partial{C_i}}{\partial{v}}\in I$.

Now, since $v$ appears effectively in $I_{C_i}$, the set
$$\mathbb{C}'=\mathbb{C}\setminus\{C_i\}\cup \left\{ \dfrac{\partial{C_i}}{\partial{v}}\right\}$$
is autoreduced and has the same rank as $\mathbb{C}$, hence $\mathbb{C}'$ is a characteristic set of $I$. Moreover,
the initial of $\dfrac{\partial{C_i}}{\partial{v}}$ is equal to $\dfrac{\partial{I_{C_i}}}{\partial{v}}$, hence
it does not depend on the leaders of $\mathbb{C}$. Yet $\dfrac{\partial{C_i}}{\partial{v}}$ is not a multiple
of $C_i$, which contradicts Corollary~\ref{Multiples}.
\end{proof}

\begin{remark} Proposition~\ref{Separants} will be extremely important for our study of the bound on the orders of the elements
of characteristic set. Section~\ref{Result} will tell about this in detail.
\end{remark}

\end{subsection}

\begin{subsection}{Another algorithm for computing canonical characteristic sets}
Let a differential ranking be fixed.

In \cite[page 636]{Unmixed} an algorithm {\sf Invert} for inverting polynomials w.r.t a characteristic set is presented. Its input consists of a polynomial $f$ and a characteristic set $\mathbb{C}$ of a characterizable ideal $(\mathbb{C}):H_\mathbb{C}^\infty$. Let $L$ and $N$ be the sets
of leaders and non-leaders of $\mathbb{C}$, respectively.
Then, the output of the algorithm is:
\begin{itemize}
\item `yes' and the inverse polynomial
$g \in k[N][L]$ such that $fg = 1 \mod (\mathbb{C})$ in the ring $k(N)[L]$, if $f$ is invertible w.r.t. $\mathbb{C}$;
\item `no', otherwise.
\end{itemize}

In the following algorithm while inverting initials we always get the answer `yes', so we just use the remaining (polynomial) part of its output. Let $\mathbb{C} = C_1,\ldots,C_p$ be an autoreduced set. For each $i,$ $1 \Le i \Le p,$ denote $\mathbb{C}(i) = C_1,\ldots,C_i$.

\begin{algorithm}\label{CanonicalAlgorithm}{\sf Canonical Characteristic Set}

{\sc Input:} a characteristic set $\mathbb{C}=C_1,\ldots,C_p$ of a characterizable differential ideal $I = [\mathbb{C}]:H_\mathbb{C}^\infty$ in
$k\{y_1,\ldots,y_l\} = k\{Y\}.$

{\sc Output:} the canonical characteristic set of $I$.

\begin{itemize}
\item Let $C_i = I_{C_i}u_{C_i}^{n_i}+a_{1,i}u_{C_i}^{n_i-1}+\ldots+a_{n_i,i},$ $1\Le i \Le p$.
\item $L := $ {\sf Leaders}$(\mathbb{C}).$
\item $N := \Theta Y\setminus \Theta L.$
\item for $i$ from $2$ to $p$ do
\begin{itemize}
\item $I_i' := $ {\sf Invert }$(I_{C_i},\mathbb{C}(i-1))$.
\item $C_i := I_i'C_i$.
\item $C_i := $ {\sf Pseudo\_Remainder}$(C_i)$ w.r.t. $\mathbb{C}(i-1).$
\end{itemize}
\item for $i$ from $1$ to $p$ do
\begin{itemize}
\item $C_i := C_i/\gcd(I_{C_i},a_{1,i},\ldots,a_{n_i,i})$ in $k[N][L].$
\item divide $C_i$ by its leading coefficient from $k$
\end{itemize}
\item Return $\mathbb{C}$.
\end{itemize}
\end{algorithm}
\end{subsection}

\begin{remark}
Note that the coefficients $a_{i,j}$ at the end of the algorithm are the {\it new} coefficients of the new $C_i$.
\end{remark}

% \begin{remark} Maple has a procedure doing the same as we do at each step of the second cycle
% of Algorithm~\ref{CanonicalAlgorithm}. It is called {\sf primpart}.
% \end{remark}

\begin{lemma}
  Algorithm~\ref{CanonicalAlgorithm} is correct, i.e., its output is the canonical characteristic set of $I$.
\end{lemma}
\begin{proof}
First, note that after each iteration of the first
{\it for}-loop the set $\mathbb{C}$ is still a characteristic set of the ideal $I$. Indeed, no $C_i$ disappears during those reductions, because the initials $I_{C_i}$ are always invertible in our situation.

The fact that the initials of $C_i$ do not depend on leaders can be proved by induction on $i$.
For $i=1$ this is the case, since $C_1$ has the lowest rank in the characteristic set. Assume that the initials of the polynomials in $\mathbb{C}(i-1)$ do not depend on leaders. Consider the initial of $C_i$, $I_i''=I_{C_i}I_i'$, after multiplication by $I_i'$.
We have $I_i''-1\in(\mathbb{C}(i-1))$ in $k(N)[L]$. During the computation of the pseudo-remainder of $C_i$ w.r.t. $\mathbb{C}(i-1)$, the
initial of $I_i''$ may be pseudoreduced by some polynomials from $\mathbb{C}(i-1)$ or multiplied by their initials, which by inductive assumption belong to $k[N]$. Hence, if $I_i'''$ is the initial of $C_i$ after the computation of the pseudo-remainder, we have
$I_i'''-f\in(\mathbb{C}(i-1))$, where $f\in k[N]$. Given that $I_i'''$ and $f$ are reduced w.r.t. $\mathbb{C}(i-1)$, we obtain
$I_i'''-f=0$, i.e., $I_i'''\in k[N]$.

The remaining two conditions for the canonicity (see Definition~\ref{SecondDefinition}) of the output of Algorithm~\ref{CanonicalAlgorithm} are ensured by the second {\it for}-loop.
\end{proof}

% \begin{remark} One can define a canonical characteristic set as the output of Algorithm~\ref{CanonicalAlgorithm}. In
% Section~\ref{DifferentDefinitions} we have already discussed another algorithm (Algorithm~\ref{Canonical}) that can also
% serve as a definition of the canonical characteristic sets. All these definitions are equivalent. This is shown in
% Proposition~\ref{Irreducible} and Theorem~\ref{HubertCanonical}.
% \end{remark}

We now know how to compute canonical characteristic sets of characterizable differential ideals from any characteristic set characterizing the ideal. The canonical characteristic set is unique and has good properties of minimality. The next section is devoted to establishing facts about
the bounds for characteristic sets. First, we do this for prime differential ideals (Theorem~\ref{OrderlyElimination}). Then, Corollary~\ref{CharacterizableOrdinaryInvertibleBound} gives us a generalization of this result to characterizable ideals. Finally, in Section~\ref{CanonicalProperties} we apply these results to canonical characteristic sets.
\end{section}

\begin{section}{Bounds for the orders of characteristic sets}\label{Result}
\begin{subsection}{Preparation}
Let $R = k\{y_1,\ldots,y_l\}$ with $\Delta = \{\delta\}$. So, we are in the ordinary case. {\it Differential dimension} of a  differential ideal $I$ is the maximal number $q$ such that $I \cap k\{y_{i_1},\ldots,y_{i_q}\} = \{0\}$.
Recall that the order of a differential polynomial $f$ is the maximal order of differential variables appearing effectively in $f$.

Fix any differential ranking. Let $\mathbb{A} = A_1,\ldots,A_p$ be an autoreduced set. Define the order of $\mathbb{A}$ by the following equality: $\ord\mathbb{A} = \ord A_1 + \ldots + \ord A_p$. Let an {\it orderly} differential ranking be fixed. If $\mathbb{C}$ is a characteristic set of a {\it prime} differential ideal $P$  then, by definition, the order of the ideal $P$ equals $\ord\mathbb{C}$.

Denote by $P(s)$ the set of elements of $P$ whose order is less than or equal to $s$. The set $P(s)$ is a prime algebraic ideal in the corresponding polynomial ring. As it is proved in \cite[II.12, Theorem 6]{Kol} or \cite[Theorems 5.4.1, 5.4.4]{Pan} the dimension of $P(s)$ is a polynomial in $s$ for $s \Ge h = \ord P$.

More precisely, $\dim P(s) = q(s+1) + \ord P$, where $q$ is the differential dimension of the ideal $P$. Moreover, $q = l-p$, where $p$ is the number of elements of a characteristic set of the ideal $P$
w.r.t. any orderly ranking. Thus, the numbers $\ord P$ and $p$ do not depend on the choice of an orderly ranking.

We are going to define the order of a characterizable differential ideal and we should be very careful because of the following example.

\begin{example} Consider the radical differential ideal $\{x(x+y')\} = I$ characterizable w.r.t. the elimination ranking $x >_{el} y$. While $I = [x]\cap[x+y']$ and the leaders of $x$ and $x+y'$ w.r.t. the ranking are the same the orders of the components are different. This is because the ideal $I$ is not characterizable w.r.t. any orderly ranking.
\end{example}

Hence, we give the following definition.

\begin{definition}\label{CharacterizableOrderDefinition} For a characterizable differential ideal $I=\bigcap\limits_{i=1}^nP_i$, where $P_i$ are minimal differential prime components of $I$, define $$\ord I = \max\limits_{1\Le i \Le n}\ord P_i.$$
\end{definition}

\begin{remark} The theory of differential dimensional polynomials is due to Kolchin \cite{Kol}. Carr\`a Ferro and Sit continued to develop this subject \cite{Sit1,CF1,CF2}. Many of the results concerning differential dimension polynomials are summarized
in \cite{Pan}. The latter book also presents many algorithms for computing these polynomials.
\end{remark}

\begin{lemma}\cite[Proposition 17]{Sadik}\label{OrderLess} Consider a prime differential ideal $P$ of differential dimension $q$ and order $h$. For every subset $\{y_{i_1},\ldots,y_{i_{q+1}}\}$ of $\{y_1,\ldots,y_l\}$, the ideal $P$ contains a differential polynomial in the indeterminates $\{y_{i_1},\ldots,y_{i_{q+1}}\}$ of order less than or equal to $h$.
\end{lemma}

A characteristic set of a prime differential ideal is not unique. For example, consider the ideal $[x] \in k\{x,y\}$ and the elimination ranking with $x > y$. Then the set $y^{(n)}x$ is a characteristic set of the ideal $[x]$ for any $n \Ge 0$. Hence, if we do not impose additional restrictions, we will not be able to obtain a bound on the order of characteristic sets of a prime differential ideal.

In order to avoid this problem, in \cite{Sadik} Sadik introduced the concept of {\it irreducible} characteristic set and proved its existence
\cite[Lemma 19]{Sadik} for any prime differential ideal.

\begin{definition}\label{SadiksIrreducible}
A characteristic set of a prime ideal is called {\it irreducible in Sadik's} sense if
\begin{enumerate}
\item $C_1$ is an irreducible polynomial in $k\{y_1,\ldots,y_l\}$,
\item each $C_{i+1}$ is irreducible in the ring $$R_i = \Quot(k[V_i]\:/\:(C_1,\ldots,C_i):I_i^\infty)[U_i],$$ where
{\begin{itemize}
\item $V_i$ is the set of all variables appearing in the polynomials $C_1,\ldots,C_i$,
\item $I_i^\infty$ is the multiplicative system generated by the initials of the polynomials $C_1,\ldots,C_i$,
\item $U_i$ is the set of variables from $C_{i+1}$ that are not in $V_i$.
\end{itemize}}
\end{enumerate}
\end{definition}

The following result allows us to find a characteristic set of a prime differential ideal with good bounds on the orders of its elements w.r.t. {\it any} differential ranking.

We formulate the following result (Lemma~\ref{IrreducibleSeparants}) in the way we are going to use it. One can conclude from its proof that this is nothing else as:
{\it a characteristic set which is irreducible in Sadik's sense
satisfies the second condition of Lemma~\ref{IrreducibleSeparants}.}

\begin{lemma}\label{IrreducibleSeparants} A prime differential ideal $P$ in $k\{y_1,\ldots,y_l\}$ admits a characteristic set $\mathbb{A} = A_1,\ldots,A_p$ with the following properties:
\begin{enumerate}
\item it is irreducible in Sadik's sense;
\item let $y_t^{(s)}$ be a differential variable appearing in $\mathbb{A}$. Assume also that $y_t^{(s)}$ does not appear in $A_1,\ldots,A_{i-1}$ but does appear in $A_i$. Then $$S_{i,t} = \dfrac{\partial A_i}{\partial y_t^{(s)}} \notin P.$$
\end{enumerate}
\end{lemma}
\begin{proof} Suppose that the second condition is failed for a characteristic set $A_1,\ldots,A_p$ irreducible in Sadik's sense,
which exists by \cite[Lemma 19]{Sadik}.
Let $z$ be a variable that does not appear in $\mathbb{A} = A_1,\ldots,A_{i-1}$ but does appear in $A_i$ and satisfies $\dfrac{\partial A_i}{\partial z}\in P$. Take the canonical characteristic set $C_1,\ldots,C_p$ of the ideal $P$.

Consider the unique factorization domain $R_{i-1}$ constructed from $A_1,\ldots,A_{i-1}$ in Definition~\ref{SadiksIrreducible}. The variable $z$ is an indeterminate in this ring. Since $\mathbb{A}$ is irreducible, the polynomial $A_i$ is irreducible in $R_{i-1}$. The polynomial $C_i$ is reducible to zero
w.r.t. $\mathbb{A}$. Hence $C_i$ is reducible to zero w.r.t. $A_i$ in $R_{i-1}$, since $A_1,\ldots,A_{i-1}$ is a characteristic set of the prime ideal $(A_1,\ldots,A_{i-1}):I_{i-1}^\infty$. Then, there exists a number $k$ and a polynomial $D_i \in R_{i-1}$ such that $I_{A_i}^kC_i = D_iA_i$. Since $D_iA_i$ is divisible by $C_i$ and $A_i$ is irreducible, we have $C_i = E_iA_i$ for some factor $E_i$ of $D_i$. Thus, the polynomial $C_i$ must contain the variable $z$. Since the polynomial $f = I_{C_i}A_i-I_{A_i}C_i\in P$
is reduced w.r.t. $A_i$, we have $$f\in
J:=(A_1,\ldots,A_{i-1}):I_{i-1}^\infty,$$ where $I_{i-1}$ denotes the multiplicative set generated by the initials $I_{A_1},\ldots,I_{A_{i-1}}.$

Since $z$ does not appear in $A_1,\ldots,A_{i-1}$, there exist generators of the ideal $J$ not containing this variable. Hence
$\dfrac{\partial f}{\partial z} \in P$. On the other hand,
$$\dfrac{\partial f}{\partial z} = \dfrac{\partial A_i}{\partial z}I_{C_i} - \dfrac{\partial {C_i}}{\partial
z}
I_{A_i} + \dfrac{\partial I_{C_i}}{\partial z}A_i - \dfrac{\partial I_{A_i}}{\partial
z}
C_i \equiv \dfrac{\partial A_i}{\partial z}I_{C_i} - \dfrac{\partial {C_i}}{\partial
z}
I_{A_i} \pmod P.$$ Thus, from $\dfrac{\partial A_i}{\partial z}\in P$ and Proposition~\ref{Separants}, we have $I_{A_i} \in P$.
But the initials of a characteristic set of a prime ideal cannot belong to it. Contradiction.
\end{proof}
\end{subsection}

\begin{subsection}{Bound for prime differential ideals}

\begin{theorem}\label{OrderlyElimination} Let $P$ be a prime differential ideal of order $h$ in $k\{y_1,\ldots,y_l\}$ and $>$ be a differential ranking ({\it not} necessarily orderly!). Then there exists a characteristic set $\mathbb{C} = C_1,\ldots,C_n$ of the ideal $P$ w.r.t. the ranking $>$ such that the order in $y_t$ of each $C_i$ does not exceed $h$ for all $1\Le t \Le l$.
\end{theorem}
\begin{proof}
For a characteristic set $\C$ of $P$ denote the set 
$$
\left\{y_k\:|\: \theta y_k\ \text{is not a leader of any}\ C_j,\ 1 \Le j \Le n, \theta \in \Theta\right\}
$$ 
by $\mathfrak{N}.$
If for some $\theta \in \Theta$ and $t,$ $1 \Le t \Le l,$ the variable $\theta y_t$ is {\it the leader} of some $C_j$ then we will show that $\ord(C_q,y_t) \Le h$
for all $1\Le q \Le n$ using Lemma~\ref{OrderLess}. Indeed, 
since $\C$ is autoreduced, we have
\begin{align}\label{ineq1}
&&\ord(C_q,y_t) \Le \ord\theta,
\end{align}
for all $q,$ $1 \Le q \Le n.$ Since $\dim P = \#\mathfrak{N},$ by Lemma~\ref{OrderLess} there exists a polynomial $$0 \ne f \in k\{y_t,\mathfrak{N}\}\cap P$$ of
order not greater than $h$. This polynomial depends only on non-leaders $\mathfrak{N}$ and the leading variable $y_t.$ Moreover, $f$ is reducible
to zero w.r.t. $\C.$ Hence, 
\begin{align}\label{ineq2}
&&\ord\theta = \ord(C_j,y_t) \Le \ord(f,y_t) \Le h.
\end{align}
The inequalities~\eqref{ineq1} and \eqref{ineq2} give us
$$
\ord(C_q,y_t) \Le h
$$
for all $q,$ $1 \Le q \Le n.$

Now let $y_t \in \mathfrak{N}$ and $\mathbb{C}$ be a characteristic set which Lemma~\ref{IrreducibleSeparants} provides to us.
Let also $y_{C_j}$ denote the differential indeterminate such that $\theta y_{C_j}$ is the leader of $C_j$ for some $\theta \in \Theta$, that is, $y_{C_j}$
is the leading variable of $C_j.$
The main idea is to reduce the polynomial {\it of the smallest order with respect to $y_{C_j}$}
$$
f_j \in k\{y_{C_j},\mathfrak{N}\}\cap P
$$ 
given by Lemma~\ref{OrderLess} w.r.t. $\mathbb{C}.$  Let $u = y_{C_j}^{(r)}$ be the
derivative of $y_{C_j}$ of the highest order in $f_j$. If we represent $f_j$
as a univariate polynomial in $u$ then denote by $I_{f_j}$ its leading coefficient. Notice that $I_{f_j}$ does not have to be the initial of
$f_j$ w.r.t. our ranking, but we still use this notation for convenience. 
For instance, $I_{f_j}$ would {\it be} the initial of $f_j$ w.r.t. the elimination
ranking $y_{C_j} > \mathfrak{N}.$
We emphasize that 
$$
I_{f_j} \notin P.
$$ 
Suppose that for some $j, 1\Le j \Le n,$ we have 
\begin{align}\label{ineq3}
&&\ord(C_j, y_t) > h.
\end{align}
%Note that 
Since $f_j$ is reducible to zero w.r.t. $\C$ we must have
\begin{align}\label{ineq4}
&&\ord\left(f_j,y_{C_j}\right) \Ge \ord\left(C_j,y_{C_j}\right).
\end{align}
Denote by ``$\arg\max\ord$''  the {\it set} of all elements which provide the maximum of the order. Consider 
$$
\tilde{\mathbb{C}} = \arg\max\limits_{C_j \in \mathbb{C}}\ord(C_j,y_t)
$$ 
and then choose $C_i \in \tilde{\mathbb{C}}$ of the {\it lowest} possible rank.
%\footnote{In \cite{Sadik} Sadik used induction that allowed him to prove the %result just for elimination rankings.}  
We can have many elements in $\tilde{\mathbb{C}}$. But we take the special one, $C_i$.
Let $u_i = \theta_i y_i$ for some $\theta_i \in \Theta$ and $u_i$ be the leader of $C_i$ for simplicity.
From \eqref{ineq3} and \eqref{ineq4} we have 
\begin{align}\label{ineq5}
&&s = \ord(C_i, y_t) > h
\end{align}
and 
$$
r_f = \ord(f_i,y_i) \Ge \ord(C_i,y_i) = r_C,
$$ where 
$$
f_i = f_i(y_i,\mathfrak{N}) = I_{f_i}\left(y_i^{(r_f)}\right)^{n_f} + a_1\left(y_i^{(r_f)}\right)^{n_f-1}+\ldots+a_{n_f}.
$$

Let us reduce each term (coefficients $a_j$, ``initial'' $I_{f_i}$ and its ``leader'' $y_i^{(r_f)}$) of $f_i$ first by $C_i$. We need to differentiate $C_i$ $q$ times and get the remainder $\tilde{f},$ where $0  \Le q \Le r_f - r_C$. Remember that $f_i$ depends only on $y_i, \mathfrak{N},$ and their derivatives.
By reduction here we mean the following. Any proper derivative $\theta$ of 
$C_i$ is linear in $\theta u_i$ and its initial is equal to the separant of
$C_i.$ We simply multiply $f_i$ by a sufficient power (say, $n_f$) of the separant and replace $y_i^{(r_f)}$ and the derivatives of $y_i$ of lower
order in $f_i$ by the corresponding tails. 

Hence, applying further steps of reduction to the terms of $\tilde{f}$ w.r.t. all $C_j$ we need to differentiate them {\it less} than $q$ times if $C_j \in \tilde{\mathbb{C}}$. Indeed, the fact that $C_i < C_j,$ as $C_i$ has the smallest rank in $\tilde\C,$ implies 
$$
\ord\left(C_i,y_{C_j}\right)<\ord\left(C_j,y_{C_j}\right).
$$
We need to differentiate them {\it not greater} than $q$ times if $C_j \notin \tilde{\mathbb{C}}$. Indeed, the set $\mathbb{C}$ is autoreduced, so
$$
\ord\left(C_i,y_{C_j}\right)\Le \ord\left(C_j,y_{C_j}\right).
$$  
In addition, the variables to reduce can come just from derivatives of variables from $C_i$.

In the case of $r_f = r_C$ we are in \cite[Lemma 20]{Sadik} because of our choice of $C_i$ and immediately get the inequality 
\begin{align}\label{hineq}
&&\ord(f_i,y_t) \Ge \ord(C_i,y_t).
\end{align}
Since $\ord(f_i,y_t) \Le h,$ the inequality~\eqref{hineq} contradicts
to inequality~\eqref{ineq5}.

Consider the other case of $r_f > r_C.$ 
Here, after we reduce all leaders of $\mathbb{C}$ from $f$ we get the polynomial depending effectively on $y_t^{(s+q)}$ and $s+q \Ge s$.
Its leading coefficient w.r.t. the variable $y_t^{(s+q)}$ is equal to 
\begin{align}\label{initial}
&&I_{C_1}^{i_1}\cdot\ldots\cdot I_{C_n}^{i_n}\cdot S_{C_1}^{j_1}\cdot\ldots\cdot S_{C_n}^{j_n}\cdot\tilde{I}_{f_i}\cdot\left(\frac{\partial C_i}{\partial y_t^{(s)}}\right)^{n_f},
\end{align}
where $i_1,\ldots,i_n,j_1,\ldots,j_n \in \mathbb{Z}_{\Ge 0}$ and $\tilde{I}_{f_i}$ is the remainder of $I_{f_i}$ w.r.t. $\mathbb{C}.$ 
Remember that $P$ is a prime ideal. Hence, 
\begin{align}\label{doesnotbelong}
&&I_{C_1}^{i_1}\cdot\ldots\cdot I_{C_n}^{i_n}\cdot S_{C_1}^{j_1}\cdot\ldots\cdot S_{C_n}^{j_n} \notin P,
\end{align}
because $I_{C_j}$ and $S_{C_j} \notin P$ for all $j,$  $1 \Le j \Le n.$
Moreover, $P = [\mathbb{C}]:H_\mathbb{C}^\infty$ and $\mathbb{C}$ is a characteristic set of $[\mathbb{C}]:H_\mathbb{C}^\infty$.
Also, 
\begin{align}\label{reduceddoesnotbelong}
&&\tilde{I}_{f_i} \notin P,
\end{align}
because $I_{f_i} \notin P$ due to our choice of $f_i.$
By the Rosenfeld lemma, the remainder of  $f_i$ we are computing
belongs to the prime algebraic ideal $(\C):H_{\C}^\infty.$
Thus, according to Lemma~\ref{LoseVar}, 
its leading coefficient given by \eqref{initial} is reducible to zero w.r.t. $\mathbb{C}.$ For a prime differential ideal the fact that an element is reducible to zero w.r.t. a characteristic set is nothing else the element belongs to the ideal. Using \eqref{doesnotbelong} and \eqref{reduceddoesnotbelong}
we conclude that the polynomial $\frac{\partial C_i}{\partial y_t^{(s)}}$
belongs to $P.$
Finally, this contradicts to Lemma~\ref{IrreducibleSeparants}.
\end{proof}

\begin{remark} Note that Theorem~\ref{OrderlyElimination} is actually a generalization of Sadik's result \cite[Theorem 24]{Sadik} that was proved just for elimination rankings. In the proof of Theorem~\ref{OrderlyElimination} we
construct the set $\tilde\C$ and choose a special element $C_i \in \tilde\C.$
Sadik used induction here. 
\end{remark}
\end{subsection}

\begin{subsection}{Characterizable ideals: main estimate for the bound}

We do not need the ordinary case for the following result.
Fix a ring of differential polynomials $k\{y_1,\ldots,y_l\}$.

\begin{theorem}\label{CharacterizableBound}
Suppose a function $h$ from the set of prime differential ideals to the set $\mathbb{Z}_{\Ge 0}$ is such that for any prime differential ideal $P$ there exists its characteristic set $C_1,\ldots,C_p$ with the property $\ord C_i \Le h(P)$ for all $i$, $1 \Le i \Le p$. Then for any characterizable differential ideal $I$ there exists its characteristic set $\mathbb{B} =  B_1,\ldots,B_k$ characterizing this ideal ($I = [\mathbb{B}]:H_\mathbb{B}^\infty$) such that
$$\ord B_i \Le \max_{1\Le j \Le n} h(P_j) =: h(I)$$
for all $i$, $1 \Le i \Le k$, where the set of ideals $\{P_j\:|\: 1 \Le j \Le n\}$ is the minimal prime decomposition of $I$.
\end{theorem}
\begin{proof}
Take the minimal prime decomposition $I = \bigcap\limits_{j=1}^n P_j$ and choose a characteristic set $\mathbb{C}_j = C_{j,1},\ldots,C_{j,p_j} \subset P_j$ with $\ord C_{j,i} \Le h(P_j) \Le h(I)$ for all $i$, $1\Le i \Le p_j,$ and $j$, $1 \Le j \Le n$.
We have
$$
I = \bigcap_{j=1}^n [\mathbb{C}_j]:H_{\mathbb{C}_j}^\infty.
$$

Let $\mathbb{B}$ be any characteristic set of $I$ characterizing this radical differential ideal, i.e., $I = [\mathbb{B}]:H_\mathbb{B}^\infty$, and $L$ be the set of its leaders which is uniquely determined by $I$ and does not depend on the choice of $\mathbb{B}$. Let $N$ be the (infinite) set of {\it all} other variables from $k\{y_1,\ldots,y_l\}$.
From Theorem~\ref{Prime} we know that
$$
J = (\mathbb{B}):H_\mathbb{B}^\infty = \bigcap_{j=1}^n (\mathbb{C}_j):H_{\mathbb{C}_j}^\infty.
$$
in the ring $k[N,L]$ and $\mathbb{B}$ is an algebraic characteristic set of $J$ which can be computed, e.g., from the reduced Gr\"obner basis $G$ of the ideal $J$. We just need to notice that $G$ can be computed from all $\mathbb{C}_j$ without involving extra variables from the set $N$. To conclude that $I = [\mathbb{B}]:H_\mathbb{B}^\infty$ we use \cite[Lemmas 3.5, 3.9, and 6.1]{Fac}.
\end{proof}

Let us switch to the {\it ordinary case} and see what corollaries Theorem~\ref{CharacterizableBound} gives us.

\begin{corollary}\label{CharacterizableOrdinaryInvertibleBound} In the ordinary case for a characterizable differential ideal $I$ there exists a characteristic set $\mathbb{C}=C_1,\ldots,C_p$ with the following properties:
\begin{itemize}
\item $I = [\mathbb{C}]:H_\mathbb{C}^\infty.$
\item $\ord C_i \Le \ord I$ (see Definition~\ref{CharacterizableOrderDefinition}) for all $i$, $1\Le i \Le p$.
\end{itemize}
\end{corollary}
\begin{proof} Follows from Theorem~\ref{OrderlyElimination} and Theorem~\ref{CharacterizableBound} setting $h(P) = \ord P$.
\end{proof}

\end{subsection}

\end{section}

\begin{section}{Bound for the order of the canonical characteristic set}
\label{CanonicalProperties}

We need the ordinary case for the following assertions about bounds. In this theorem we also reprove the uniqueness of the canonical characteristic set of a characterizable differential ideal. Note that we need the ordinary case only for our bound but not for the uniqueness. Let some differential ranking be fixed.

\begin{theorem}\label{HubertCanonical} Let $\mathbb{C} = C_1,\ldots,C_p$ be the canonical characteristic set of a
characterizable differential ideal $I$. The set $\mathbb{C}$ is uniquely determined and has the following bound $$\ord C_i \Le \ord I$$
for all $i,1\Le i \Le p.$
\end{theorem}
\begin{proof} Let $\mathbb{B}$ be any characteristic set of the ideal $I$ with invertible initials and separants. Such a set characterizes $I$, i.e., we have $I = [\mathbb{B}]:H_{\mathbb{B}}^\infty$. Hence, $$[\mathbb{C}]:H_{\mathbb{C}}^\infty = [\mathbb{B}]:H_{\mathbb{B}}^\infty.$$

Let us move into the ring $k(N)[L]$. Both $\mathbb{C}$ and $\mathbb{B}$ are autoreduced
sets of the same rank. Hence, they have the same leaders and the degrees of these leaders. Thus, all the elements of $\mathbb{B}$ are partially reduced w.r.t. $\mathbb{C}$.
Then, $\mathbb{B} \subset (\mathbb{C}):H_{\mathbb{C}}^\infty = (\mathbb{C})$ by
Theorem~\ref{Rosenfeld}. The latter equality follows from \cite[Proposition 3.3 and Lemma 6.1]{Fac}

The initials and separants of $\mathbb{B}$ are invertible. Hence, for any $h \in H_\mathbb{B}^\infty$ there exists $h' \in k[N][L]$ such that $hh' = 1 + f,$ where $f \in (\mathbb{B})$. Let $a \in (\mathbb{B}):H_{\mathbb{B}}^\infty$. So, there exists $h \in H_\mathbb{B}^\infty$ such that $ha \in (\mathbb{B})$. Hence, $hh'a = (1+f)a=a+af \in (\mathbb{B})$ and $a \in (\mathbb{B})$.

Then, $(\mathbb{B}):H_{\mathbb{B}}^\infty = (\mathbb{B}) \subset (\mathbb{C}):H_{\mathbb{C}}^\infty$. In the same way, $(\mathbb{C}):H_{\mathbb{C}}^\infty  \subset (\mathbb{B}):H_{\mathbb{B}}^\infty.$ Thus,
$$(\mathbb{C}):H_{\mathbb{C}}^\infty  = (\mathbb{B}):H_{\mathbb{B}}^\infty.
$$
As a result, Algorithm~\ref{Canonical} does give its answer
independently of the choice of the input characteristic set of a
fixed characterizable differential ideal.

Let us prove now the bound for the orders. Compute the canonical characteristic set of the ideal $I$ using the result of Corollary~\ref{CharacterizableOrdinaryInvertibleBound}. So, let $\mathbb{B} = B_1,\dots,B_p$ be a characteristic set given by Corollary~\ref{CharacterizableOrdinaryInvertibleBound}. We have $\ord B_i \Le \ord P$ for all $i,1\Le i \Le p.$

No steps of Algorithm~\ref{Canonical} involve differentiations and no new variable can appear in the output. Since the canonical characteristic set does not depend on the choice of input characteristic set of $I$, we have the bound for $\mathbb{C}$ we need.
\end{proof}

\end{section}

\begin{section}{Computation of the canonical characteristic set from generators}

We do {\it not} assume the ordinary case now.  Fix a differential ranking. The algorithm for computing a characteristic set of a prime differential ideal (represented by its generators as a radical differential ideal, $I = \{F\}$) is given in \cite[Theorem 6]{Bou1}. It simply
\begin{itemize}
\item takes the first component of its {\sf Rosenfeld\_Gr\"obner} decomposition (appeared in \cite{Bou1});
\item computes reduced Gr\"obner basis of the correspondent algebraic ideal;
\item applies algebraic pseudo-autoreduction to this basis, extracting a characteristic set of the ideal $I$.
\end{itemize}

Note that the above algorithm works only for prime differential ideals. The case of {\it characterizable} ideals specified by {\it sets of generators}
as radical differential ideals is more tricky. The following Algorithm~\ref{CharCharComputation} computes the canonical characteristic set in this case (see Proposition~\ref{FromGenerators} for the proof).

\begin{remark}
It is not known how to perform the inverse transformation between the two representations, i.e., how to compute a set of generators of a characterizable ideal given its canonical characteristic set.
\end{remark}

\begin{algorithm}\label{CharCharComputation}{\sf Characteristic Set of a Characterizable Differential Ideal}

{\sc Input:} a finite set $F$ of differential polynomials generating a characterizable differential ideal.

{\sc Output:} the canonical characteristic set of $\{F\}$.

\begin{itemize}
\item Let $\mathfrak{C} =$ {\sf Rosenfeld\_Gr\"obner}$(F)$ and $\mathfrak{C} = \mathbb{C}_1,\ldots,\mathbb{C}_n$.
\item Let $[\mathbb{C}_{i_j}]:H_{\mathbb{C}_{i_j}}^\infty$ be the components whose characteristic sets have the sets of leaders
 of the highest possible rank in $\mathfrak{C}$ and $1 \Le j \Le k$.
\item Let $I' = \bigcap\limits_{j=1}^k\left(\mathbb{C}_{i_j}\right):H_{\mathbb{C}_{i_j}}^\infty$.
\item $L := $ {\sf Leaders}$\left(\mathbb{C}_{i_1}\right).$
\item $N := \Theta Y\setminus \Theta L.$
\item $GB := $ {\sf Reduced\_Gr\"obner\_Basis}$(I')$ in $k(N)[L]$.
%\item $N' := \{x\in N\:|\: x\ \mathrm{appears\ in}\ GB\}.$
\item $\mathbb{D} := $ {\sf Clear\_out\_denominators}$(GB)$ in $k(N)[L].$
\item divide each element of $\mathbb{D}$ by its leading coefficient from $k$.
\item Return $\mathbb{D}$.
\end{itemize}
\end{algorithm}

One can think that the property of a radical differential ideal to be characterizable can be checked in this way (by computing the canonical characteristic set). But this is not the case. Radical differential ideals having characteristic sets satisfying Definition~\ref{SecondDefinition} may not be characterizable. This is illustrated in Examples~\ref{DiffExample} and \ref{AlgExample}.

\begin{remark} Algorithm {\sf Rosenfeld\_Gr\"obner} is presented in \cite{Bou1} and \cite{Bou2} and implemented in Maple.
\end{remark}

\begin{remark} Note that in the second line of the above algorithm it would not be sufficient to consider only the characterizable
components having characteristic sets of the highest rank in $\mathfrak{C}$. Indeed, let $x>y>z$, and consider
the following algebraic characterizable ideal and its decomposition into characterizable components:
$$I=(y^2+z, x^3+x^2y+xy-z) = (y^2+z, x+y) \cap (y^2+z,x^2+y).$$
The characteristic sets of both components have the same set of leaders, $\{x,y\}$. The component
of the highest rank is $(y^2+z,x^2+y)$ and, clearly, $I\neq (y^2+z,x^2+y)$.
\end{remark}

\begin{proposition}\label{FromGenerators} Algorithm~\ref{CharCharComputation} computes the canonical characteristic set of the given characterizable differential
ideal $\{F\}$.
\end{proposition}
\begin{proof}
Let $\mathbb{C}$ be the canonical characteristic set of the characterizable ideal $I=\{F\}$.

First, let us prove an auxiliary
\begin{lemma}\label{Inclusion}
  Let $P$ be a prime differential ideal with a characteristic set $\mathbb{A}$ whose set of leaders coincides with that
  of $\mathbb{C}$, where $\mathbb{C}$ is a characteristic set of $[\mathbb{C}]:H_{\mathbb{C}}^\infty = I$. Assume also that $I\subseteq P$. Then
  $(\mathbb{C}):H_{\mathbb{C}}^\infty\subset(\mathbb{A}):H_{\mathbb{A}}^\infty$.
\end{lemma}
\begin{proof}
   Let $f\in(\mathbb{C}):H_{\mathbb{C}}^\infty$.
   Then $f$ is partially reduced w.r.t. $\mathbb{C}$. Since the leaders
   of $\mathbb{A}$ and $\mathbb{C}$ coincide, $f$ is partially reduced w.r.t. $\mathbb{A}$.
   Since $f\in I$ and $I\subseteq P$, we have $f\in P$.
   Hence, by Rosenfeld's lemma (Theorem~\ref{Rosenfeld}), $f\in(\mathbb{A}):H_{\mathbb{A}}^\infty$.
\end{proof}

Consider the prime decomposition $I=\bigcap P_i$, where $P_i$'s are the essential prime components
of $I$. Let $\mathbb{A}_i$ be a characteristic set of $P_i$, then, according to Theorem~\ref{Prime}, ideal $P_i'=(\mathbb{A}_i):H_{\mathbb{A}_i}^\infty$
is a minimal prime component of the algebraic ideal $(\mathbb{C}):H_{\mathbb{C}}^\infty$.

Consider also the essential prime decompositions $J_l=\bigcap Q_{lj}$ of the characteristic components
$J_l=[\mathbb{C}_l]:H_{\mathbb{C}_l}^\infty$ of $I$. The intersection of these decompositions is a finite prime decomposition
of $I$. According to \cite[Section I.16]{Rit}, every essential prime component appears in every
finite prime decomposition of the radical ideal $I$, which implies that every $P_i$ can be found among $Q_{lj}$.
Moreover, according to Theorem~\ref{T20}, the leaders of $\mathbb{A}_i$
coincide with the leaders of $\mathbb{C}$, hence $P_i$ can be found among those
$Q_{lj}$ whose characteristic sets have leaders coinciding with the leaders of $\mathbb{C}$.

Applying Theorem~\ref{T20} again, we obtain that $P_i$ can be found among the
essential prime components of those $J_l$ whose characteristic sets $\mathbb{C}_l$ have leaders coinciding with
the leaders of $\mathbb{C}$. Now, since for each $l$, $I\subseteq J_l$, the rank of the
set of leaders of $\mathbb{C}_l$ is lower than or equal to the rank of the set of leaders of $\mathbb{C}$.
Hence, $P_i$ can be found among the essential prime components of those $J_l$, for which
the set of leaders of $\mathbb{C}_l$ has the highest rank, i.e., among the essential prime
components of $J_{i_1},\ldots,J_{i_k}$.

Thus, by Theorem~\ref{Prime}, every
minimal prime $P_i'$ of the algebraic ideal  $(\mathbb{C}):H_{\mathbb{C}}^\infty$
can be found among the minimal primes of the algebraic ideals
$(\mathbb{C}_{i_1}):H_{\mathbb{C}_{i_1}}^\infty,\ldots,(\mathbb{C}_{i_k}):H_{\mathbb{C}_{i_k}}^\infty$,
and we obtain $$(\mathbb{C}):H_{\mathbb{C}}^\infty\supseteq\bigcap_{j=1}^k (\mathbb{C}_{i_j}):H_{\mathbb{C}_{i_j}}^\infty=I'.$$
The inverse inclusion follows from the above Lemma~\ref{Inclusion}.
Hence, $I'=(\mathbb{C}):H_{\mathbb{C}}^\infty$, and the canonical characteristic set $\mathbb{D}$ of
$I'$ computed by the above algorithm coincides with that of $(\mathbb{C}):H_{\mathbb{C}}^\infty$ and
of $I$.
\end{proof}

\end{section}

\begin{section}{Canonical characteristic sets of prime differential ideals}

Prime differential ideals are characterizable. Can we say more about their canonical characteristic sets?
Theorem~\ref{OrderlyElimination} is based on finding a characteristic set of a prime differential ideal satisfying two properties:
\begin{enumerate}
\item Separants are not in the ideal\label{first};
\item Sadik's property of irreducibility\label{second}.
\end{enumerate}

Such a characteristic set exists (see Lemma~\ref{IrreducibleSeparants}). Moreover, the orders of its elements are bounded by the order of the ideal. Applying Algorithm~\ref{CanonicalAlgorithm}, we get a canonical characteristic set certainly satisfying the same bound on the orders.

So, it is natural to ask whether this uniquely determined canonical characteristic set satisfies the properties~\ref{first}--\ref{second}. It turns out that the answer is only partially positive! The first property holds for a bigger class of ideals, namely, characterizable ones. This is shown in Proposition~\ref{Separants}. The second property is not necessarily true. This is shown in Example~\ref{Counterexample}.

\begin{example}\label{Counterexample} Consider the ideal $I = \{x^2-t,(zx+1)y+1\} \subset k\{x,t,z,y\}$ and any differential ranking such that $y>z>x>t$. The set $x^2-t,(zx+1)y+1$ is a characteristic set of $I$ and it is irreducible in Sadik's sense.

Nevertheless, the canonical characteristic set of $I$, which is equal to $x^2-t,(z^2t-1)y+zx-1$ is not irreducible by Sadik because
\begin{align*}
(z^2t-1)y+zx-1 &= (z^2x^2-1)y+zx-1 = \\
&= (zx-1)(zx+1)y+(zx-1)=\\
&= (zx-1)((zx+1)y+1)\\
\end{align*}
 in the polynomial ring $\Quot(k[x,t]/(x^2-t))[y,z].$
\end{example}

\begin{remark} Example~\ref{Counterexample} also shows that Proposition~\ref{LowestRank} is not true if we omit the condition that the initials depend only on non-leaders.
\end{remark}

In conclusion, we note that the object we computed, i.e. the canonical characteristic set of a prime differential ideal, is not just unique but also has additional natural properties such as a low differential order of its elements, as well as other properties showing that the choice of requirements for the canonical characteristic set is not arbitrary. The importance of a characteristic set having these properties is, again, shown in Theorem~\ref{OrderlyElimination}.

\end{section}

\begin{section}{Testing equality of characterizable ideals} \label{Equality}

As was mentioned in the Introduction, the canonical characteristic set
trivially allows to check equality of two characterizable ideals: two
characterizable ideals are equal if and only if their canonical
characteristic sets coincide. However, one can also check the equality
of characterizable ideals specified by any characterizing sets using
the following criterion.

\begin{theorem}\label{EqualityCriterion}
  Let $\mathbb{A}$ and $\mathbb{B}$ be characteristic sets of
  characterizable ideals $I=[\mathbb{A}]:H_\mathbb{A}^\infty$ and 
  $J=[\mathbb{B}]:H_\mathbb{B}^\infty$, respectively.
  Then $I=J$ iff $\mathbb{A}\subset J$ and $\mathbb{B}\subset I$.
\end{theorem}
\begin{proof}
Assume that $\mathbb{A}\subset J$, $\mathbb{B}\subset I$. 
Since $\mathbb A$ is a characteristic set of $I$, i.e., an
autoreduced subset of $I$ of the least rank, and since
$\mathbb B$ is an autoreduced subset of $I$, we have
$\rank{\mathbb A}\le\rank{\mathbb B}$. Symmetrically,
$\rank{\mathbb B}\le\rank{\mathbb A}$. Thus, the ranks
of $\mathbb A$ and $\mathbb B$ are equal. 

%Show first
%that $$\rank_\le \mathbb{A}=\rank_\le \mathbb{B}.$$ Let $\mathbb{A} = A_1,\ldots,A_q$
%and $\mathbb{B}=B_1,\ldots,B_p$. If $\rank_\le A_1 < \rank_\le B_1$ then $A_1$ is
%not reducible w.r.t. $\mathbb{B}$ that contradicts to the inclusion $\mathbb{A}\subset J$
%and to the fact that $\mathbb{B}$ characterizes the ideal $J$. The
%similar consideration gives $\rank_\le A_1 = \rank_\le B_1$. Assume
%that $\rank_\le A_1,\ldots,A_i = \rank_\le B_1,\ldots B_i$. If 
%$\rank_\le A_{i+1} < \rank_\le B_{i+1}$ then $A_{i+1}$ is not
%reducible w.r.t. $\mathbb{B}$, because of the induction hypothesis and
%the fact that $\mathbb{A}$ is autoreduced. Thus, again $\rank_\le
%A_{i+1} = \rank_\le B_{i+1}$ and finally $\rank_\le \mathbb{A}=\rank_\le \mathbb{B}$.

We show now that $I=J$. It is sufficient to demonstrate that $I\subset
J$, since the proof of the inverse inclusion is symmetric. Moreover, 
since $J$ is characterized by $\mathbb{B}$, it is sufficient to show
that for all $f\in I$, $f$ can be reduced to $0$ w.r.t. $\mathbb{B}$. 
Let $f\in I$, and let $f'$ be the differential remainder of $f$
w.r.t. $\mathbb{B}$. Then, since $\mathbb{B}\subset I$, $f'\in I$. But
$f'$ is not reducible w.r.t. $\mathbb{B}$ and, given that 
$\rank_\le \mathbb{A}=\rank_\le \mathbb{B}$,
$f'$ is also not reducible w.r.t. $\mathbb{A}$. Hence, $f'=0$.

The converse implication, that $I=J$ implies inclusions
$\mathbb{A}\subset J$ and $\mathbb{B}\subset I$, is straightforward.
\end{proof}

\begin{proposition}\label{HNotIn} If $\mathbb{A}$ is a characteristic set of a
non-trivial characterizable differential ideal $[\mathbb{A}]:H_\mathbb{A}^\infty = I$
then $H_\mathbb{A} \notin I$.
\end{proposition}
\begin{proof}
Assume that $H_\mathbb{A}\in I$. Then there exists $h\in
H_\mathbb{A}^\infty$
such that $h\cdot H_\mathbb{A}\in [\mathbb{A}]$. Hence,
for $h_1=h\cdot H_{\mathbb{A}}\subset H_{\mathbb A}^\infty$ we have $1\cdot h_1\in [\mathbb{A}]$,
which implies that $I$ contains $1$. 
\end{proof}

%Assume that $H_\mathbb{A} \in I$. Then, $H_\mathbb{A}$ belongs to each prime component of $I$. Hence,
%one of the initials or separants $H$ of
%$\mathbb{A}$ is in some minimal prime component $P$ of $I$. But by Theorem~\ref{KandriInvertibleSep} this $H$ must be invertible, since $\mathbb{A}$ is a characteristic set of the
%characterizable differential ideal $I$. So, there exist polynomials $f$ and  $g$ with the latter one depending only on non-leaders such that $H\cdot f \equiv g \pmod I$. We have $g\in P$. According to Theorems~\ref{T20}~and~\ref{Prime} any characteristic set of $P$ has the same set of leaders as one of $I$. Thus, $g$ is not reducible w.r.t. a characteristic set of $P$. Contradiction.
%\end{proof}

\begin{corollary}
The criterion given by Kolchin:
``two prime differential ideals $I=[\mathbb{A}]:H_\mathbb{A}^\infty$
and $J = [\mathbb{B}]:H_\mathbb{B}^\infty$ given by their
characteristic sets $\mathbb{A}$
and $\mathbb{B}$ coincide iff
\begin{enumerate}
\item $\mathbb{A} \subset J$,
\item $\mathbb{B} \subset I$,
\item $H_\mathbb{A} \notin J$,
\item $H_\mathbb{B} \notin I$;''
\end{enumerate}
holds in the case of
characterizable differential ideals.
\end{corollary}
\begin{proof}
Necessity can be derived from Proposition~\ref{HNotIn}. Sufficiency follows from Theorem~\ref{EqualityCriterion}.
\end{proof}

\end{section}

\begin{section}{Examples}\label{Examples}

We illustrate how Algorithm~\ref{CanonicalAlgorithm} works.

\begin{example} Consider again the ideal $I = \{x^2-t,(zx+1)y+1\} \subset k\{x,t,z,y\}$ and any differential ranking such that $y>z>x>t$. The set $x^2-t,(zx+1)y+1$ is a characteristic set of $I$. According to algorithm {\sf Invert }  we introduce a new variable $w$ and compute the reduced Gr\"obner basis $GB$ of the ideal $(x^2-t,w-(zx+1))$ w.r.t. the lexicographic ordering with $x>w>z>t$.

We have $$GB = w^2-2w-z^2t+1, zx+1-w, xw-x-zt, x^2-t$$ and
$$GB\cap k[w,z,t] = w^2-2w-z^2t+1 = P(w).$$ Substituting $w=0$ we get $P(0) = 1-z^2t \ne 0.$ Hence, $zx+1$ is invertible and $$(zx+1)(zx+1-2) \equiv (zx+1)^2-2(zx+1) \equiv z^2t-1 \mod (x^2-t).$$
Take $$(zx+1-2)((zx+1)y+1) \equiv (z^2t-1)y + zx-1 \mod (x^2-t).$$
So, the canonical characteristic set is equal to
$$
x^2-t, (z^2t-1)y + zx-1.
$$
\end{example}

We know that any characterizable ideal has the canonical characteristic set. Is the converse true? If a radical ideal $I$ has a characteristic set satisfying the conditions of Definition~\ref{SecondDefinition} this does not imply that the ideal $I$ is characterizable! If we involve differentiations that the most common example in constructive differential algebra comes into the play:

\begin{example}\label{DiffExample} Consider the ring of differential polynomials $k\{y\}$ and the radical differential ideal
$$I := \{y'^2+y\} = [y'^2+y]:(y')^\infty\cap [y],$$
which is not characterizable, but its characteristic set $y'^2+y$ satisfies all the conditions of  Definition~\ref{SecondDefinition}.
\end{example}

In the case of zero dimensional radical algebraic ideals the converse is true. More precisely, from Lemma~\cite[Lemma 3.5]{Fac} it follows that if a zero dimensional radical ideal has the canonical characteristic set then this ideal is characterizable. This is not true in non-zero dimensional cases.

\begin{example}\label{AlgExample}
Consider the polynomial ring $k[x,y,z,t]$, the ranking $y > x > z > t,$ and the radical ideal
$$
I = (zx+t, zy + t) = (zx+t,zy+t):z^\infty\cap(z,t) = (zx+t,x-y)\cap(z,t).
$$
One cannot exclude the second component $(z,t)$ from this decomposition. So, it is minimal and the ideal $I$ is not characterizable by \cite[Theorem 3]{Ovch} (see also \cite{Ovch2}). Nevertheless, its characteristic set $zx+t, zy + t$ w.r.t. the ranking $y > x > z > t$ is canonical.
\end{example}

%We return to the discussion of the method of testing equality from \cite[Exercise 1, page 171]{Kol}.
%
%\begin{example}\label{CharNotPrime} Consider the set $\mathbb{C} = x(x-1),y,z$ and the  ideal $(\mathbb{C}) = I$ characterizable 
%w.r.t. the ranking %$z>y>x$. Let $A_1=x(x-1), A_2=xy, A_3=(x-1)z$. The set $\mathbb{A} = A_1,A_2,A_3$ is a characteristic 
%set of the same ideal $I$ but the product of the %initials $x\cdot(x-1) \in I$. Thus, the criterion does not work.
%\end{example}

\end{section}

\begin{section}{Conclusions}
The summary of results from \cite{Bou2,Bou3,Fac,Pol,Dif} and
the new results obtained in this paper suggest that
the canonical characteristic set yields
a convenient representation of a characterizable
differential ideal: its non-algorithmic definition is simple
and transparent, it allows to solve algorithmic problems such
as ideal membership, it satisfies several natural properties,
including the {\it bound on the orders} of its elements,
and there exist efficient algorithms that convert other representations
of a characterizable ideal (by a regular system or by any other
characterizing set) into the canonical one. Moreover, one can
compute the canonical characteristic set of a characterizable 
differential ideal, if the latter is specified either by its 
{\it set of generators} as a radical differential ideal, or 
by its (possibly redundant) characteristic decomposition.  
\end{section}

\begin{acknowledgements}
We are highly grateful to Michael F. Singer for suggesting to us
the directions of exploring the subject of this paper. 
We thank Evelyne Hubert for the idea that has lead to the present 
definition and proof of uniqueness of a 
canonical differential characteristic set, for
useful references, and for a productive discussion.
We are thankful to Evgeniy V. Pankratiev for helpful 
comments and support, to William Y. Sit and 
Alexey Zobnin for important discussions, and also thank 
Alexey Zobnin for extensive comments on the draft of this paper.
We appreciate suggestions of the referees very much.
\end{acknowledgements}

\bibliographystyle{amsplain}
\bibliography{canonical}

\providecommand{\bysame}{\leavevmode\hbox to3em{\hrulefill}\thinspace}
\providecommand{\MR}{\relax\ifhmode\unskip\space\fi MR }
\providecommand{\MRhref}[2]{%
  \href{http://www.ams.org/mathscinet-getitem?mr=#1}{#2}
}
\providecommand{\href}[2]{#2}
\begin{thebibliography}{10}

\bibitem{Wei}
T.~Becker and V.~Weispfenning, \emph{Gr{\"o}bner bases}, Springer-Verlag, New
  York-Berlin-Heidelberg, 1993.

\bibitem{Bou1}
F.~Boulier, D.~Lazard, F.~Ollivier, and M.~Petitot, \emph{Representation for
  the radical of a finitely generated differential ideal}, Proceedings of ISSAC
  1995, ACM Press, 1995, pp.~158--166.

\bibitem{Bou2}
\bysame, \emph{Computing representations for radicals of finitely generated
  differential ideals}, Tech. report, IT-306, LIFL, 1997.

\bibitem{Bou3}
F.~Boulier and F.~Lemaire, \emph{Computing canonical representatives of regular
  differential ideals}, Proceedings of ISSAC 2000, ACM Press, 2000, pp.~38--47.

\bibitem{Bou4}
F.~Boulier, F.~Lemaire, and M.~Moreno-Maza, \emph{{PARDI!}}, Proceedings of
  ISSAC 2001, ACM Press, 2001, pp.~38--47.

\bibitem{Unmixed}
D.~Bouziane, A.~Kandri~Rodi, and H.~Ma\^arouf, \emph{Unmixed-dimensional
  decomposition of a finitely generated perfect differential ideal}, Journal of
  Symbolic Computation \textbf{31} (2001), 631--649.

\bibitem{Bb}
B.~Buchberger, \emph{An algorithm for finding a basis for the residue class
  ring of a zero-dimensional polynomial ideal (in german)}, Ph.D. thesis,
  University of Insbruck, 1965.

\bibitem{CF1}
G.~Carr{\`a}~Ferro, \emph{Some properties of the lattice points and their
  application to differential algebra}, Communications in Algebra \textbf{15}
  (1987), no.~12, 2625--2632.

\bibitem{CF2}
\bysame, \emph{Some remarks on the differential dimension}, Lecture Notes in
  Computer Science \textbf{357} (1989), 152--163.

\bibitem{ComputingResolvent}
T.~Cluzeau and E.~Hubert, \emph{Computing resolvent representation of a regular
  differential ideal}, In preparation, 2006.

\bibitem{DahanPoster}
X.~Dahan, M.~Moreno~Maza, E.~Schost, and X.~Jin, \emph{Change of ordering for
  regular chains in positive dimension}, Poster Presentation, ISSAC, 2006.

\bibitem{Bounds}
O.~Golubitsky, M.~Kondratieva, M.~Moreno~Maza, and A.~Ovchinnikov, \emph{Bounds
  for algorithms in differential algebra}, In preparation, 2006.

\bibitem{Fac}
E.~Hubert, \emph{Factorization-free decomposition algorithms in differential
  algebra}, Journal of Symbolic Computation \textbf{29} (2000), no.~4-5,
  641--662.

\bibitem{Pol}
\bysame, \emph{Notes on triangular sets and triangulation-decomposition
  algorithms {I}: {P}olynomial systems.}, Symbolic and Numerical Scientific
  Computing 2001, 2003, pp.~1--39.

\bibitem{Dif}
\bysame, \emph{Notes on triangular sets and triangulation-decomposition
  algorithms {II}: {D}ifferential systems}, Symbolic and Numerical Scientific
  Computing 2001, 2003, pp.~40--87.

\bibitem{Hubert04}
\bysame, \emph{Improvements to a triangulation-decomposition algorithm for
  ordinary differential systems in higher degree cases}, Proceedings of ISSAC
  2004, ACM Press, 2004, pp.~191--198.

\bibitem{Kol}
E.R. Kolchin, \emph{Differential algebra and algebraic groups}, Academic Press,
  New York, 1973.

\bibitem{Pan}
M.V. Kondratieva, A.B. Levin, A.V. Mikhalev, and E.V. Pankratiev,
  \emph{Differential and difference dimension polynomials}, Kluwer Academic
  Publisher, 1999.

\bibitem{Ovch}
A.~Ovchinnikov, \emph{Characterizable radical differential ideals and some
  properties of characteristic sets}, Programming and Computer Software
  \textbf{30} (2004), no.~3, 141--149.

\bibitem{Ovch2}
\bysame, \emph{On characterizable ideals and characteristic sets},
  Contributions to General Algebra \textbf{14} (2004), 91--108.

\bibitem{Rit}
J.F. Ritt, \emph{Differential algebra}, American Mathematical Society, New
  York, 1950.

\bibitem{Sadik}
B.~Sadik, \emph{A bound for the order of characteristic set elements of an
  ordinary prime differential ideal and some applications}, Applicable Algebra
  in Engineering, Communication and Computing \textbf{10} (2000), no.~3,
  251--268.

\bibitem{Sit1}
W.~Y. Sit, \emph{Differential dimension polynomials of finitely generated
  extensions}, Proceedings of American Mathematical Society \textbf{68} (1978),
  no.~3, 251--257.

\bibitem{Sit}
W.Y. Sit, \emph{The {R}itt-{K}olchin theory for differential polynomials},
  Differential Algebra and Related Topics, Proceedings of the International
  Workshop (NJSU, 2--3 November 2000), 2002.

\end{thebibliography}

\end{document}